

\baselineskip=14pt
\parskip=10pt

\font\eightrm=cmr8 
\font\eighttt=cmtt8
\magnification=\magstephalf

\def\1{{\overline{1}}}
\def\2{{\overline{2}}}
\parindent=0pt
\overfullrule=0in

\def\frac#1#2{{#1 \over #2}}
\bf
\centerline
{
Computational and Theoretical Challenges On Counting Solid Standard Young Tableaux
}
\rm
\bigskip
\centerline{ {\it
Shalosh B. EKHAD and
Doron 
ZEILBERGER}\footnote{$^1$}
{\eightrm  \raggedright
Department of Mathematics, Rutgers University (New Brunswick),
Hill Center-Busch Campus, 110 Frelinghuysen Rd., Piscataway,
NJ 08854-8019, USA.
{\eighttt zeilberg  at math dot rutgers dot edu} ,
\hfill \break
{\eighttt http://www.math.rutgers.edu/\~{}zeilberg/} .
First written: Feb. 20, 2012. This version: Feb. 28, 2012 [thanks to Manuel Kauers and Fredrik Johansson] .
Exclusively published in {\eighttt http://www.math.rutgers.edu/\~{}zeilberg/pj.html} and {\tt arxiv.org} .
Supported in part by the USA National Science Foundation.
}
}

{\bf Very Important:} This article is accompanied by the Maple package \hfill\break
{\eighttt http://www.math.rutgers.edu/\~{}zeilberg/tokhniot/SolidSYT} \quad .
Readers who have access to Maple should download it and read it into a Maple session,
so that they can follow the text more vividly.

The subject of {\it solid partitions} goes back to Percy MacMahon [M], but to our surprise, 
as far as we can tell by {\it googling}, no one has seriously studied {\it Solid Standard Young Tableaux}.
Let's first recall some basic facts about the familiar kind.

{\bf Review of (2D) Standard Young Tableaux}

Recall that a (usual) (2D) {\it Young diagram} of shape $\lambda=(\lambda_1, \dots, \lambda_k)$ 
(where $\lambda_1 \geq \lambda_2 \geq \dots \geq \lambda_k >0$
are integers)
is a left-justified collection of $k$ rows of empty unit-boxes, where the top row has
$\lambda_1$ boxes, the second row has $\lambda_2$ boxes, $\dots$, and the bottom, $k$-th row, has
$\lambda_k$ boxes. For example, the following is a Young diagram of shape $(3,2,2,1)$
$$
\matrix{
X&X&X\cr
X&X\cr
X&X\cr
X\cr} \quad ,
$$
where $X$ denotes an empty unit-box. Let $n:=\lambda_1+ \dots + \lambda_k$ be the number of boxes (alias the
``integer that is being partitioned'' by $\lambda$). A {\it Standard Young Tableau} is a way of
placing the integers $1$ through $n$ inside the boxes, so that 
{\bf all rows} and {\bf all columns} are
{\bf increasing} (when read from left-to-right, and top-to-bottom respectively). For example,
$$
\matrix{
1&3&8\cr
2&5\cr
4&6\cr
7\cr} \quad ,
$$
is one of the seventy Standard Young Tableaux of Shape $(3,2,2,1)$. There is a beautiful formula, due
to Frame, Robinson, and Thrall[FRT], called the {\it hook-length formula}, for the number of
Standard Young Tableaux of a given shape.  Calling that number $f_{\lambda}$, it is:
$$
f_{\lambda}=
\frac{n!}{\prod_{b \in \lambda} h_b} \quad ,
$$
where for any {\it box} $b$
in the Young diagram, $h_b$ is its {\it hook-length}, i.e. 
the number of boxes that
are either weakly to its right or weakly below. For 
example, the hook-lengths for the shape $(3,2,2,1)$ are
$$
\matrix{
6&4&1\cr
4&2\cr
3&1\cr
1\cr} \quad .
$$
It follows, that
$$
f_{3221}=
\frac{8!}{6 \cdot 4 \cdot 1 \cdot 4 \cdot 2 \cdot 3 \cdot 1 \cdot 1} = 70
\quad ,
$$
as claimed above.

The {\bf proof from the book}
of that amazing formula is due to Curtis Greene, Albert Nijenhuis, and Herbert Wilf[GNW]. 
In addition to its considerable {\it face-value}, this proof is also historically significant,
since it turned DZ from an unhappy analyst into a happy combinatorialist!

The [GNW] {\it proof-from-the-book} is presumbly a serendipitous {\it by-product} of
the Greene-Nijenhuis-Wilf {\it algorithm-from-the-book} (also told in [GNW])
to generate a Standard Young Tableau of a given shape {\it uniformly at random}. 
It goes like this.

First roll an $n$-faced {\it fair} die, and decide accordingly the starting box. Then whenever visiting
a box $b$, roll a fair $(h_b-1)$-faced die, and decide which box in the hook (except the one you
are at right now, you must move on!) to go to next. Keep doing it, until you wind-up at a corner, where
there is nowhere to go. Put $n$ there. Now you have a smaller shape, with $n-1$ empty boxes, and
continue recursively.

[This algorithm is implemented in procedure {\tt GNW} of the Maple package  \hfill\break
{\eighttt http://www.math.rutgers.edu/\~{}zeilberg/tokhniot/GreeneNijenhuisWilf} ] .

The total number of Standard Young Tableaux with $n$ cells, $\omega_n:=\sum_{\lambda \dashv n} f_{\lambda}$
is the famous ``number of involutions'', sequence {\tt http://oeis.org/A000085} \quad ,
(thanks to the so-called {\it Robinson-Schenstead algorithm}), that has a very simple,  {\it recurrence}
$\omega_n=\omega_{n-1}+(n-1)\omega_{n-2}$, that enables one to easily compute the first ten thousand terms in a fraction of a second.

{\bf Solid Standard Young Tableaux}

Now the shapes are the 3D Young diagrams of {\it plane partitions}. 
Recall that a {\it plane partition} is a 
two-dimensional array of positive integers $p_{ij}$ where both
rows and columns are {\it weakly-decreasing} 
and its {Young diagram} consists of piling $p_{ij}$ empty boxes
above location $(i,j)$ on the floor. 
A  {\it Solid Standard Young Tableau} of a given shape (with $n$ empty boxes)
is a way of placing the integers $1$ through $n$ such that 
going from left-to-right, from back-to-front, and
from down-to-up the entries are {\bf increasing}.

It is unrealistic to expect a {\it nice} formula
for the number of Solid 
Standard Young Tableaux of a given 3D shape, but we can still easily compute it, using
an obvious recurrence. If the number of boxes in our shape is $n$, then the entry $n$ must reside in one of
the {\it corners} (unit-boxes where all the forward-going neighbors (in each of the three directions) do not belong to the shape).
Then $f_{\lambda}$ (where $\lambda$ is now a 3D shape (alias a plane-partition)) is the sum of
$f_{\lambda'}$  over all $\lambda'$ obtained by removing a corner box from $\lambda$. This is implemented in 
{\tt SolidSYT}'s procedure {\tt Nu(L)}, where $L$ is a shape (i.e. plane partition), expressed as a list-of-lists
of positive integers.

The 3D analog of $\omega_n$, obtained by 
summing $Nu(L)$ over all plane partitions $L$ of a given integer $n$,
is implemented by procedure {\tt Y3number(n);}, and to get the
first $K$ terms of that sequence, type {\tt Y3numberSeq(K);}.

For the record, here are the first thirty terms, 
[taken from 
{\eighttt http://www.math.rutgers.edu/\~{}zeilberg/tokhniot/oSolidSYT1}]:
$1, 3, 9, 33, 135, 633, 3207, 17589, 102627, 636033, 4161141, 28680717, 207318273, 1567344549, 12345147705,$
$101013795753, 856212871761,7501911705747, 67815650852235, 631574151445665, 6051983918989833,$
$    59605200185016639, 602764245172225251, 6252962956009863363,   66482211459036254169, 723810526382641418667, $
$8062440364611311185977, 91804267420894431624357, 1067720130017504052805449,    12673922788286515247094267 \quad .$

{\bf Generating a Uniformly-at-Random Solid Standard Young Tableau of a Given Shape}

Using the beautiful approach of [W], clearly explained and exploited in [NW], one can generate
{\it uniformly at random}, a Solid Standard Young Tableau of a given (solid) shape.
If $\lambda$ is such a shape (alias plane partition) with $n$ boxes, then the entry $n$ can reside in any of its
{\it corners} (boxes where none of its forward-going neighbors are in the shape). Let the set
of corners be $C$. Then we form a {\bf loaded die} whose faces are labeled by the members of $C$,
and the probability of it lending on face $c$ is $f_{\lambda - c}/f_{\lambda}$, where 
$f_\lambda$ is the number of Solid Standard Young Tableaux of shape $\lambda$ (implemented by procedure
{\tt Nu(L)} in {\tt SolidSYT}). We then place the $n$ in the corner-box decided by the die, and
get a smaller shape, $\lambda'=\lambda - c$ with $n-1$  boxes, and continue recursively, until we get
the empty shape.

[Procedure {\tt RSSYT(L);} implements this in the Maple package {\tt SolidSYT}, try for example, \hfill\break
{\tt RSSYT([[3,3,3],[3,3,3],[3,3,3]]);} for getting, {\it uniformly-at-random},
one of the  $6405442434150$ ways
of placing 1 through 27 in a $3 \times 3 \times 3$
box, in such a way that when you go form left-to-right, from 
back-to-front, and from down-to-up, they
are always increasing.]

{\bf The Three-Dimensional Greene-Nijenhuis-Wilf Algorithm}

As the shapes get larger, {\tt RSSYT} gets slower and slower, since it relies on the recursive procedure
{\tt Nu} (there is no (known) closed-form expression for $f_\lambda$ for three-dimensional shapes).

The beauty of the Greene-Nijenhuis-Wilf algorithm is that it is so much faster! The die cast
at every step is always fair! Unfortunately, the three-dimensional analog no longer gives you
a random Solid Standard Young Tableau {\it uniformly}, but is gets (experimentally) fairly close.
So if you don't mind a little bias, you are welcome to use procedure {\tt GNW3(L);} .

For example

{\tt GNW3([[10\$10]\$10]);} gives you, instantaneously, a (not-quite-uniformly-at) random way of placing 
$1$ through $1000$ in a $10 \times 10 \times 10$
box, in such a way that when you go form left-to-right, from back-to-front, and from down-to-up, they
are always increasing.

{\bf Some Computational Challenges Regarding the Enumeration of Solid Standard Young-Tableaux of Cylindrical Shapes}

Most of us know that the number of (usual, 2D) Standard Young Tableaux of shape $(n,n)$ is given by
the famous {\it Catalan  Numbers} $(2n)!/(n!(n+1)!)$, {\tt http://oeis.org/A000108} .
The number of Standard Young Tableaux of shape $(n,n,n)$ is given by the so-called
{\it three-dimensional Catalan Numbers}, {\tt http://oeis.org/A005789} \quad , that count the number of
ways of walking $n$ steps in a 3D Manhattan always staying in $x \geq y \geq z$,
and (as of Feb. 18, 2012) Sloane has it up to the five-dimensional version {\tt http://oeis.org/A005791} \quad .

Since (what Sloane calls) the $k$-dimensional Catalan numbers are
given by the explicit (hypergeometric!) formula
$(k-1)!(nk)!/(n! \cdots (n+k-1)!)$ (that follows immediately from
the hook length formula)
it follows that
the enumerating sequence $a(n)$ (for each specific $k$)
is a {\it hypergeometric}
sequence, in other words, there exist polynomials $p_1(n)$ and $p_0(n)$ such that
$$
p_0(n)a(n)+p_1(n)a(n+1)=0 \quad,
$$
which is a special case (first-order) of a very important {\it ansatz}, the so-called {\it holonomic}, or $P$-{\it recursive}
ansatz, that consists of sequences satisfying a linear-recurrence equation with polynomial coefficients of some (finite) order $L$:
$$
\sum_{i=0}^{L} p_i(n)a(n+i)=0 \quad .
$$

Going to sequences enumerating Solid Standard Young Tableaux of Cylindrical Shapes, i.e. $\lambda \times \{1, \dots , n \}$,
for a (usual) partition $\lambda$, things are very mysterious. For $\lambda=(2,1)$, i.e. for 3D-shapes of
the form $[[n,n],[n]]$ we have (not-quite-so-trivially!) the famous {\it Kreweras sequence},
{\tt http://oeis.org/A006335} (why?), that is also a hypergeometric sequence, and hence holonomic.

Procedure {\tt Sidra(L,n,N0);} in our Maple package {\tt SolidSYT} spits out the first {\tt N0} terms of
the enumerating sequence for a family of shapes $L$ with parameter $n$. For example, to get the first
20 terms of the Kreweras sequence type: \quad  {\tt Sidra([[n,n],[n]],n,20); }  \quad .

On the other hand, if $\lambda=[[n,n],[n,1]]$,  
the sequence ``{\tt Sidra([[n,n],[n,1]],n,40); }''
can be ``described'' (empirically, so far) by a
second-order linear recurrence equation with polynomial coefficients, see:

{\tt http://www.math.rutgers.edu/\~{}zeilberg/tokhniot/oSolidSYT3} \quad  .

{\bf 1st Rigorous Challenge} ($0.01$ US dollars) : Find a rigorous proof of this recurrence .

To our surprise, the enumerating sequence for the number of Solid Standard Young Tableaux of the cylindrical shapes
$(2,1,1) \times \{1, \dots n\}$, i.e., in the notation of {\tt SolidSYT}, 

{\tt Sidra([[n,n],[n],[n]], n, BigEnough);}

for which, with some effort, (with {\tt BigEnough}=120), we were able to find the first $121$ terms, see

{\tt http://www.math.rutgers.edu/\~{}zeilberg/tokhniot/oSolidSYT7} \quad ,

did not yield a linear recurrence equation with polynomial coefficients of order {\tt ORDER} and degree {\tt DEGREE} with
{\tt (ORDER+1)(DEGREE+1)} less than 115. This brings us to the:

{\bf 1st Non-Rigorous Challenge} ($100$ US dollars) : Find a linear recurrence equation with polynomial coefficients
(empirically) satisfied by the sequence $a(n):=$ ``number of Solid Standard Young Tableaux'' of shape 
$(2,1,1) \times \{1, \dots n\}$ (or equivalently, $(3,1) \times \{1, \dots n\}$).
Equivalently, $a(n)$ is the number of ways of walking from $(0,0,0,0)$ to $(n,n,n,n)$
using positive  unit steps in the four-dimensional Manhattan lattice, in such a way
that all the visited points $(x_1,x_2,x_3,x_4)$ always satisfy $x_1 \geq x_2 \geq x_3$ {\bf and} $x_1 \geq x_4$.

{\bf 2nd Rigorous Challenge} ($1$ US dollar) :  Having ``conjectured'' the above recurrence
(i.e. proven it experimentally), find a ``rigorous'' proof.

{\bf Feb. 28, 2012 Update}: Manuel Kauers and his student, Fredrik Johansson kindly informed me
that no recurrence exists with {\tt (ORDER+1)(DEGREE+1)} less than 3000. See Manuel Kauers's message in: \hfill\break
{\tt http://www.math.rutgers.edu/\~{}zeilberg/mamarim/mamarimhtml/ssytAppendix.html} \quad ,

If the answer is negative, then we have the alternative

{\bf 2'nd Rigorous Challenge} ($10$ US dollars) :  Prove that the above sequence is {\bf not} holonomic. Finally, the 

{\bf Big Question}: (Lots of glory but no cash) Characterize all partitions $\lambda$ for which the sequence 
enumerating Solid Standard Young Tableaux of shape $\lambda \times \{1, \dots , n \}$, ($n=1,2, \dots $)
satisfy a linear recurrence equation with polynomial coefficients.

{\bf References}

[FRT] J.S. Frame, G. de B. Robinson, and R.M. Thrall,
{\it The hook graphs of the symmetric group}, Canad. J. Math. {\bf 6} (1954), 317-324.

[GNW] Curtis Greene, Albert Nijenhuis, and Herbert Wilf,
{\it A probabilistic proof of a formula for the
number of Young Tableaux of a given shape},
Advances in Mathematics {\bf 31}(1979), 104-109.

[M] P.A. MacMahon,  {\it ``Combinatory Analysis''}, Vol. 2, Cambridge University Press, London and New York,
1916 (reprinted by  Chelsea, New York, 1960).

[NW] Albert  Nijenhuis and Herbert Wilf,
{\it ``Combinatorial Algorithms''},
2nd edition, Academic Press,  1978.

[W] Herbert Wilf,
{\it A unified setting for sequencing, ranking, and selection
algorithms for combinatorial objects},
Advances in Mathematics {\bf 24}(1977), 281-291.

\end